\documentclass[english]{amsart}

\usepackage{babel}
\usepackage{amstext}
\usepackage{amsmath}
\usepackage{amsfonts}
\usepackage{amssymb}
\usepackage{latexsym}
\usepackage{ifthen}
\usepackage[all,2cell,dvips]{xy} \UseAllTwocells \SilentMatrices
\usepackage{xy}
\xyoption{all}
\newtheorem{lemma1}[equation]{}

\newenvironment{lemma}{\begin{lemma1}{\bf Lemma.}}{\end{lemma1}}

\newenvironment{theorem}{\begin{lemma1}{\bf Theorem.}}{\end{lemma1}}
\newenvironment{proposition}{\begin{lemma1}{\bf Proposition.}}{\end{lemma1}}
\newenvironment{corollary}{\begin{lemma1}{\bf Corollary.}}{\end{lemma1}}

\newenvironment{notation}{\begin{lemma1}{\bf Notation.}}{\end{lemma1}}

\newenvironment{conjecture}{\begin {lemma1}{\bf Conjecture.}}{\end{lemma1}}

\newcommand{\C}{\ensuremath{\mathbb{C}}}

\newcommand{\PP}{\ensuremath{\mathbb{P}}}
\newcommand{\D}{\ensuremath{\mathbb{D}}}

\newcommand{\Homsheaf} { \ensuremath{ \mathcal{H} \! om}}
\newcommand{\Ssheaf} { \ensuremath{ \mathcal{S} } }

\newcommand{\Osheaf} { \ensuremath{ \mathcal{O} } }

\newcommand{\chow}[1]{\ensuremath{\mathcal{C}(#1)}}

\newcommand{\merom}[3]{\ensuremath{#1:#2 \dashrightarrow #3}}

\newcommand{\holom}[3]{\ensuremath{#1:#2  \rightarrow #3}}
\newcommand{\fibre}[2]{\ensuremath{#1^{-1} (#2)}}

\makeatletter
\ifnum\@ptsize=0  \fi \ifnum\@ptsize=2  \fi 

\newcommand\sO{{\mathcal O}}

\newcommand\Hom{{\rm Hom}}

\DeclareMathOperator*{\rk}{rk}

\newcommand{\Chow}[1]{\ensuremath{\mathcal{C}(#1)}}

\title {The structure of uniruled manifolds with split tangent bundle} 
\date{21st December 2006}

\begin{document}

\begin{abstract}
In this paper we show that a uniruled manifold with a split tangent bundle admits almost
holomorphic fibrations that are related to the splitting. We analyse these fibrations in detail 
in several special cases, this yields new results about the integrability of
the direct factors and the universal covering of the manifold.
\end{abstract}

\author{Andreas H\"oring}
\address{Andreas H\"oring, c/o Lehrstuhl Mathematik I, Mathematisches Institut,
Universit\" at Bayreuth,
D-95440 Bayreuth, Germany}
\email{andreas.hoering@uni-bayreuth.de}

\maketitle

\tableofcontents

\section{Introduction}

A compact K\"ahler manifold $X$ has a split a tangent bundle if 
$T_X = V_1 \oplus V_2$, where $V_1$ and $V_2$ are subbundles of $T_X$.
Initiated by Beauville's conjecture \ref{conjecturebeauville} 
on the universal covering of these manifolds \cite{Bea00}, 
these manifolds have been studied by several authors during the last years (\cite{Dru00}, \cite{CP02}, \cite{BPT04},
\cite{H05}). One of the main themes of these papers is that uniruled manifolds with split tangent bundle
play a distinguished role. For example if $X$ is projective and not uniruled, 
then both $V_1$ and $V_2$ are integrable
\cite[Thm. 1.3]{H05}, while for uniruled manifolds it is easy to construct examples where this is not the case.

The goal of this paper is to develop a structure theory for 
uniruled K\"ahler manifolds of arbitrary dimension. 
We will observe in proposition \ref{propositionrationallyconnectedungeneric} 
that if $Z$ is a general fibre of the rational quotient map of $X$,
then
\[
T_Z = (T_Z \cap V_1|_Z) \oplus (T_Z \cap V_2|_Z).
\]
In particular $Z$ is a rationally connected manifold with a (maybe trivial) splitting of the
tangent bundle. This ``ungeneric position property'' (cf. \cite{H06} for the terminology)
puts us in a much better situation since we have the following description for rationally connected manifolds 
with split tangent bundle.

\begin{theorem} \cite[Thm. 1.4]{H05}
\label{theorembogomolovbetter}
Let $X$ be a rationally connected manifold such that $T_X = V_1 \oplus V_2$. 
If $V_1$ or $V_2$ is integrable,  
$X$ is isomorphic to a product $X_1 \times X_2$ such that
$V_j = p_{X_j}^* T_{X_j}$ for $j=1,2$. 
In particular both $V_1$ and $V_2$ are integrable.
\end{theorem}

So far there are no  
examples of rationally connected
manifolds with split tangent bundle where the direct factors are not integrable. 
In fact I am fairly optimistic that such examples do not exist.

\begin{conjecture}
\label{conjecturerationallyconnectedintegrable}
Let $X$ be a projective manifold with split tangent bundle $T_X=V_1 \oplus V_2$.
If $X$ is rationally connected, $V_1$ or $V_2$ is integrable.
\end{conjecture}
 
Using theorem \ref{theorembogomolovbetter} we can 
show the existence of a meromorphic fibration on $X$ that is related to the 
decomposition of the tangent bundle. More precisely we have the

\begin{theorem}
\label{theoremfibration}
Let $X$ be a uniruled compact K\"ahler manifold such that $T_X = V_1 \oplus V_2$.
Let $Z$ be a general fibre of the rational quotient map, and suppose
that $T_Z \cap V_1|_Z$ or $T_Z \cap V_2|_Z$ is integrable.
Then for $i=1,2$ there exists an almost holomorphic fibration
\merom{\phi_i}{X}{Y_i} such that the general fibre $F_i$ is rationally connected and
\[
T_{F_i} =  (T_{Z} \cap V_i|_Z)|_{F_i} \subset V_i|_{F_i}.
\]
\end{theorem}

We expect that the statement still holds without the hypothesis on the integrability,
this would of course follow from conjecture \ref{conjecturerationallyconnectedintegrable}.
We therefore show the conjecture in section \ref{sectionrationalquotient} for a splitting in small rank,
cf. lemma \ref{lemmauniruledintegrability}.
If we specify to the case where one of the direct factors
has rank 2, we obtain a more precise statement.

\begin{theorem}
\label{theoremcorank2}
Let $X$ be a uniruled compact K\"ahler manifold such that $T_X = V_1 \oplus V_2$ and $\rk V_1=2$.
Let $Z$ be a general fibre of the rational quotient map, and suppose
that $T_Z \cap V_1|_Z$ or $T_Z \cap V_2|_Z$ is integrable.
Then there are three possibilities: 
\begin{enumerate}
\item $T_Z \cap V_1|_Z = V_1|_Z$. Then the manifold 
$X$ admits the structure of an analytic fibre bundle $X \rightarrow Y$ such that the general fibre
is rationally connected and $T_{X/Y}=V_1$.
\item $T_Z \cap V_1|_Z$ is a line bundle. 
Then there exists an equidimensional map $\holom{\phi}{X}{Y}$ 
such that the general $\phi$-fibre $F$ is a rational curve and $T_F \subset V_1|_F$. 
\item $T_Z \subset V_2|_Z$.
\end{enumerate} 
\end{theorem}

One of the goals of this structure theory is to
``contract the obstruction to being integrable'', that is to construct a 
fibration $X \rightarrow Y$ such that $Y$ and the general fibre $F$ have a split tangent bundle 
with integrable direct factors. 
We realize this goal in the projective case for a splitting in vector bundles of small rank. 

\begin{theorem}
\label{theoremrank2}
Let $X$ be a uniruled projective manifold such that $T_X= \oplus_{j=1}^k V_j$, where
for all $j \in \{ 1, \ldots, k \}$ we have $\rk V_j \leq 2$. 
Let $Z$ be a general fibre of the rational quotient map. 
If $T_Z \cap V_j|_Z \neq 0$, the direct factor $V_j$ is integrable.

Furthermore the rational quotient map can be realised
as a flat fibration $\holom{\phi}{X}{Y}$ on a projective manifold
such that
\[
T_Y = \oplus_{j=1}^k (\phi_* (T \phi (V_j)))^{**}.
\]
In particular $(\phi_* (T \phi (V_j)))^{**}$ is an integrable subbundle of $T_Y$ for every $j \in \{ 1, \ldots, k \}$
(cf. \ref{notationimage} for the notation).
\end{theorem}

We then come back to the origin of our study of manifolds with split tangent bundle which is the

\begin{conjecture} (A. Beauville)
\label{conjecturebeauville}
Let $X$ be a compact K\"ahler manifold such that $T_X = V_1 \oplus V_2$, where $V_1$ and $V_2$ are vector bundles.
Let \holom{\mu}{\tilde{X}}{X} be the universal covering of $X$. Then $\tilde{X} \simeq X_1 \times X_2$, where
$p_{X_j}^* T_{X_j} \simeq \mu^* V_j$. If moreover $V_j$ is integrable, then 
there exists an automorphism  of $\tilde{X}$
such that we have an identity of subbundles of the tangent bundle $\mu^* V_j = p_{X_j}^* T_{X_j}$. 
\end{conjecture}

This will be done in section \ref{sectionuniversalcovering} where we obtain the

\begin{theorem} 
\label{theoremcorank2projective}
Let $X$ be a uniruled projective manifold such that $T_X = V_1 \oplus V_2$ and $\rk V_1=2$.
Let $Z$ be a general fibre of the rational quotient map, then one of the following holds.
\begin{enumerate}
\item $T_Z \cap V_1|_Z \neq 0$. If $V_1$ and $V_2$ are integrable, conjecture \ref{conjecturebeauville} holds.
\item $T_Z \cap V_1|_Z = 0$. Then $\det V_1^*$ is pseudoeffective and $V_2$ is integrable.
\end{enumerate} 
\end{theorem}

{\bf Acknowledgements.}
I want to thank L. Bonavero and T. Peternell for encouraging me to work on this problem. 
I also want to thank the DFG-Schwerpunkt \lq Globale Methoden in der Komplexen
Geometrie\rq \ and the \lq Deutsch-franz\"osische Hochschule-Universit\'e franco-allemande\rq \ for
partial financial support.

\section{Notation and basic results}
\label{sectionbasics}

We work over the complex field $\C$. For standard definitions in complex algebraic geometry 
we refer to  \cite{Ha77} or \cite{Kau83}, 
for positivity notions of vector bundles we follow the definitions from \cite{Laz04b}. 
Manifolds and varieties are always supposed
to be irreducible.
A fibration is a proper surjective morphism \holom{\phi}{X}{Y} with connected fibres
from a complex manifold to a normal complex variety $Y$ 
such that $\dim X > \dim Y$. The $\phi$-smooth locus is the largest Zariski open subset $Y^* \subset Y$
such that for every $y \in Y^*$, the fibre $\fibre{\phi}{y}$ is a smooth variety of dimension $\dim X - \dim Y$.
The $\phi$-singular locus is its complement.
A fibre is always a fibre in the scheme-theoretic sense, a set-theoretic fibre is the reduction of the fibre.

A meromorphic map \merom{\phi}{X}{Y} from a compact K\"ahler manifold to a normal K\"ahler variety
is an almost holomorphic fibration if  
there exist non-empty open subsets $X^* \subset X$ and $Y^* \subset Y$ 
such that \holom{\phi|_{X^*}}{X^*}{Y^*} is a fibration.
In particular for $y \in Y$ a general point, the fibre \fibre{\phi}{y} exists in the usual sense
and is compact.

Let $X$ be a compact K\"ahler manifold such that $T_X=V_1 \oplus V_2$. Suppose that $X$
is the blow-up \holom{\mu}{X}{X'} of a compact K\"ahler manifold $X'$ along a smooth submanifold.
Since in the complement of the exceptional locus we have an isomorphism $\mu^* \Omega_{X'} \simeq \Omega_X$, 
we can consider the reflexive sheaves $W_i:=(\mu_* V_i)^{**}$ as subsheaves of $T_{X'}$. Since the image
of the exceptional locus has codimension at least 2, we obtain a splitting
\[
T_{X'} = W_1 \oplus W_2.
\]
Furthermore we have an easy lemma relating the universal coverings of $X$ and $X'$.

\begin{lemma}
\label{lemmasplittingblowup}
Let $X$ be a compact K\"ahler manifold such that $T_X=V_1 \oplus V_2$. Suppose that $X$
is the blow-up \holom{\mu}{X}{X'} of a compact K\"ahler manifold $X'$ along a smooth submanifold.
Then we have a splitting $T_{X'} = W_1 \oplus W_2$.
If $W_1$ and $W_2$ are integrable and conjecture \ref{conjecturebeauville} holds for $X'$,
then the conjecture holds for $X$.
\end{lemma}

{\bf Proof.}  The statement can be shown following the proof of \cite[Prop. 4.24]{H05} and we refrain
from repeating the lengthy argument. $\square$

\begin{notation}
\label{notationimage}
Let \holom{\phi}{X}{Y} be a fibration between complex manifolds.
The canonical map $\phi^* \Omega_Y \rightarrow \Omega_X$ induces a 
a generically surjective sheaf homomorphism 
\holom{T \phi}{T_X}{\phi^* T_Y}. In particular for $\Ssheaf \subset T_X$ a quasicoherent subsheaf,
we have an inclusion $T \phi(\Ssheaf) \subset \phi^* T_Y$. Since $\phi$ is proper,
we can push forward to obtain a quasicoherent subsheaf $\phi_* (T \phi (\Ssheaf)) \subset T_Y$. 
\end{notation}

\begin{lemma}
\label{lemmapushdown}
Let $X$ be a complex manifold such that $T_X = V_1 \oplus V_2$.
Let \holom{\phi}{X}{Y} be a fibration 
onto a complex manifold $Y$ that makes $X$ into a $\PP^1-$ or conic bundle. 
Then for $j=1,2$, the reflexive sheaf $W_j:= (\phi_* (T \phi(V_j)))^{**} \subset T_Y$ 
is a subbundle of $T_Y$ and 
\[
T_Y=W_1 \oplus W_2.
\]
\end{lemma}

{\bf Proof.} 
If $\phi$
is a $\PP^1$-bundle the morphism is smooth, so \cite[Lemma 4.22]{H05} applies.
If $\phi$ is a conic bundle it is well-known
that the set $D \subset Y$ such that for $y \in D$, the fibre $\fibre{\phi}{y}$ is not reduced,
has codimension at least 2  \cite[Prop. 1.8.5]{Sar82}. Therefore \cite[Lemma 4.22]{H05} applies again. $\square$

\medskip

We recall some basic statements about holomorphic foliations, for more details we refer
to \cite{CLN85,H06}.
Let $X$ be a compact K\"ahler manifold. 
A subbundle $V \subset T_X$ is integrable if it is closed under the Lie bracket.
We recall that the Lie bracket
\[
[.,.] \ : \ V  \times V \ \rightarrow \ T_X
\]
is a bilinear antisymmetric mapping that is not $\Osheaf_X$-linear but 
induces an $\Osheaf_X$-linear map $\wedge^2 V \rightarrow T_X / V$ that is zero if and only if $V$ is integrable.
In particular 
\[
H^0(X, \Homsheaf( \wedge^2 V, T_X / V)) = 0
\] 
implies that $V$ is integrable. In general we will show this vanishing property using a dominating family $S$ 
(i.e. through a general point of $X$ passes at least one member of the family) of
subvarieties $(Z_s)_{s \in S}$ of $X$ such that a general member of the family satisfies
\[
H^0(Z_s, \Homsheaf( \wedge^2 V, T_X / V)|_{Z_s}) = H^0(Z_s, ((\wedge^2 V)^* \otimes (T_X / V))|_{Z_s}) =  0.
\] 
Since an antiample vector bundle does not have any global sections, we will use
this frequently in the following form.

\begin{lemma}
\label{lemmaintegrability}
Let $X$ be a compact K\"ahler manifold, and let $V \subset T_X$ be a subbundle.
Let $(Z_s)_{s \in S}$ be a dominating family of $X$ such that for 
a general member $Z_s$ 
of the family, the restriction of the vector bundle $\wedge^2 V|_{Z_s}$ is ample and $(T_X/V)|_{Z_s}$
is trivial. Then $V$ is integrable. $\square$
\end{lemma} 

By the Frobenius theorem an integrable subbundle $V$ of $T_X$ induces a foliation on $X$, 
i.e. for every $x \in X$ there exists
an analytic neighbourhood $U$ and a submersion $U \rightarrow W$ such that $T_{U/W}=V|_U$. 
These submersions are called the distinguished maps of the foliation and the fibres are the so-called
plaques.
This shows that $V$ can be realised
as the tangent bundle of locally closed subsets of $X$. The foliation induces an equivalence relation
on $X$, two points being equivalent if and only if they can be connected by chains of smooth (open) curves $C_i$
such that $T_{C_i} \subset V|_{C_i}$. An equivalence class is  called a leaf of the foliation. 
A subset of $X$ is $V$-saturated if it is a union of leaves.

The next proposition, which is a corollary of the global stability theorem for foliations on K\"ahler
manifold (cf. \cite{Pe01} for a short proof) 
gives a first idea why rationally connected manifolds are so useful in this context.

\begin{proposition} 
\label{propositionquotientsubmersion}
Let $X$ be a compact K\"ahler manifold such that $T_X=V_1 \oplus V_2$. 
Suppose that $V_1 \subset T_X$ is integrable and that one leaf
is compact and rationally connected. Then $X$ has the structure of an analytic fibre bundle
$X \rightarrow Y$ over a compact K\"ahler manifold $Y$ such that $T_{X/Y}=V_1$.
\end{proposition}

{\bf Proof.} By \cite[Cor. 2.11]{H05} there exists a submersion 
$X \rightarrow Y$ onto a compact K\"ahler manifold $Y$ such that $T_{X/Y}=V_1$ and the fibres
are rationally connected.
The arguments of \cite[Lemma 3.19]{H05} 
(which do not use the projectiveness hypothesis made there) 
then establish that the submersion is locally trivial.
$\square$

\section{The rational quotient map}
\label{sectionrationalquotient}

In this section we show proposition \ref{propositionrationallyconnectedungeneric} which
is the crucial observation of this paper. The moral idea behind the statement is that the
rational quotient map reflects the existence of an ample subsheaf $\Ssheaf$ of the tangent bundle $T_X$.
Proposition \ref{propositionrationallyconnectedungeneric} can then be seen as a translation of the
basic fact that a direct sum of sheaves is ample if and only if both direct factors are ample.
Once we have established this technical statement, we can use theorem \ref{theorembogomolovbetter} to show
theorem \ref{theoremfibration} and with some extra effort theorem \ref{theoremcorank2}.

\begin{proposition} 
\label{propositionrationallyconnectedungeneric}
Let $X$ be a compact K\"ahler manifold such that $T_X = V_1 \oplus V_2$.
Let $X \dashrightarrow Y$ be an almost holomorphic fibration 
such that the general fibre is rationally connected.
Then the general fibre $Z$ satisfies
\[
T_Z = (T_Z \cap V_1|_Z) \oplus (T_Z \cap V_2|_Z). 
\]
\end{proposition}

{\bf Proof.}
Let \holom{q_1}{T_X|_Z}{V_1|_Z} and \holom{q_2}{T_X|_Z}{V_2|_Z} be the restrictions of
the natural projections to $Z$. Let \holom{i_1}{V_1|_Z}{T_X|_Z}
and  \holom{i_2}{V_2|_Z}{T_X|_Z} be the inclusion maps. Furthermore we have
a natural inclusion \holom{i}{T_Z}{T_X|_Z} and a surjective map of vector bundles
\holom{q}{T_X|_Z}{N_{Z/X}}.
Taking the composition of these maps we obtain the following {\it in almost all cases not
commutative} diagram.
\[
\xymatrix{
& & V_1|_Z \ar[d] \ar[rd]^\beta & &
\\
0 \ar[r] & T_Z  \ar[r]^{i} \ar[ru]^\alpha \ar[rd]_\gamma & T_X|_Z \ar[r]^q
\ar[u] \ar[d] & N_{Z/X} \ar[r] & 0
\\
& & V_2|_Z \ar[u] \ar[ru]_\delta & &
}
\]
Note that for every $x \in Z$ we have 
$\ker \alpha_x = T_{Z,x} \cap V_{2,x}$ and $\ker \gamma_x= T_{Z,x} \cap V_{1,x}$, 
since $(T_{Z,x} \cap V_{2,x}) \cap (T_{Z,x} \cap V_{1,x})=0$ this implies
\[
(*) \rk \alpha_x + \rk \gamma_x = \dim T_{Z,x} - \dim (T_{Z,x} \cap V_{2,x}) + \dim T_{Z,x} 
- \dim (T_{Z,x} \cap V_{1,x}) \geq \rk T_{Z,x},
\]
in particular $\rk \alpha + \rk \gamma \geq \rk T_Z$.
We claim that in fact $\rk \alpha + \rk \gamma = \rk T_Z$. Granting this for the time being, we show how this 
allows to conclude.
Since $\ker \alpha = T_{Z} \cap V_{2}$ and $\ker \gamma= T_{Z} \cap V_{1}$,
we see that $\rk (T_{Z} \cap V_{1}) + \rk (T_{Z} \cap V_{2}) = \rk T_Z$
and clearly $(T_{Z} \cap V_{1}) \cap (T_{Z} \cap V_{2})=0$. Thus we have a splitting of $T_Z$ at least in
a general point.
Yet by $(*)$
we have $\rk \alpha_x + \rk \gamma_x \geq T_{Z,x}$ for every $x \in Z$, so the lower semicontinuity
of the rank of a morphism of locally free sheaves shows that $\alpha$ and $\gamma$ are of constant rank.
It follows that $T_Z \cap V_1|_Z$ and $T_Z \cap V_2|_Z$ are subbundles of $T_Z$ such that
\[
T_Z = (T_Z \cap V_1|_Z) \oplus (T_Z \cap V_2|_Z). 
\]

{\it Proof of the claim.}
The general fibre has a trivial normal bundle and is rationally connected, so 
\[
\Hom(T_Z, N_{Z/X}) \simeq H^0(Z, \Omega_Z \otimes \sO_Z^{\oplus \dim X - \dim Z})=0.
\]
Hence the morphisms \holom{\beta \circ \alpha}{T_Z}{N_{Z/X}}
and \holom{\delta \circ \gamma}{T_Z}{N_{Z/X}} are zero, in particular
\[
\rk \alpha \leq \rk (\ker \beta) = \rk V_1 - \rk \beta
\]
and
\[
\rk \gamma \leq \rk (\ker \delta) = \rk V_2 - \rk \delta.
\]
The morphism $\beta$ is the restriction of $q$ to $V_1|_Z$ and $\delta$ is the restriction of $q$ to $V_2|_Z$.
Since $q$ is surjective this implies that 
\[
\rk \beta + \rk \delta \geq \rk N_{Z/X} = \rk T_X|_Z - \rk T_Z.
\]
Putting these inequalities together we obtain
\begin{eqnarray*}
\rk T_Z & \leq & \rk \alpha + \rk \gamma
\\
& \leq & \rk V_1 + \rk V_2 - \rk \beta - \rk \delta
\\
& \leq & \rk T_X - (\rk T_X- \rk T_Z) = \rk T_Z.
\end{eqnarray*}
This proves the claim.
$\square$

{\bf Remark.}
The reader will have noticed that the proof does not really use the rational connectedness of $Z$,
but merely the cohomological condition $h^0(Z, \Homsheaf(T_Z, N_{Z/X})) = h^0(Z, \Omega_Z)=0$. 
In fact the proposition is part of a more ``ungeneric position'' theory describing fibre spaces with
split tangent bundle that is developed in \cite{H06}.
The same cohomological condition 
was used in \cite[4.4.]{Bea00} to show a more special result.

\medskip

{\bf Proof of theorem \ref{theoremfibration}.}
Let $Z$ be a general fibre of the rational quotient map of $X$. By proposition 
\ref{propositionrationallyconnectedungeneric} we have
\[
T_Z = (T_Z \cap V_1|_Z) \oplus (T_Z \cap V_2|_Z).
\]
By hypothesis one of the intersections $T_Z \cap V_1|_Z$ or $T_Z \cap V_2|_Z$
is integrable. Therefore by theorem \ref{theorembogomolovbetter} the general fibre $Z$ is isomorphic
to a product $Z_1 \times Z_2$ such that $T_Z \cap V_j|_Z=p_{Z_j}^* T_{Z_j}$ for $j=1,2$.

If $T_Z \cap V_1=0$ the identity map $X \rightarrow X$ satisfies the statement, so 
we suppose without loss of generality that $T_Z \cap V_1$ is not zero.
Since this holds for a general fibre, the 
submanifolds $Z_1 \times z$ for $z \in Z_2$ form a dominant family of submanifolds of $X$.
Let $Y_1^* \subset \Chow{X}$ be the open subset 
parametrizing the family in the cycle space,
let $\Gamma_1 \subset Y_1^* \times X$ be the graph of the family, and let \holom{q_1}{\Gamma_1}{Y_1^*} and
\holom{p_1}{\Gamma_1}{X} be the natural projections. By construction $p_1$ is dominant and an isomorphism
on its image $p_1(\Gamma_1)$, 
so we have a holomorphic map \merom{\phi_1^*:=q_1 \circ p_1^{-1}}{p_1(\Gamma_1)}{Y_1^*}.
Let $Y_1$ be the normalisation of the closure of $Y_1^* \subset \chow{X}$, then 
we obtain the stated almost holomorphic fibration \merom{\phi_1}{X}{Y_1}.
The general fibre $F_1$ of this map is just a member of the family $Z_1 \times z$, so clearly 
$T_{F_1} \subset V_1|_{F_1}$ and $F_1$ is rationally connected. The statement for $T_Z \cap V_2$ follows
analogously. $\square$

\medskip

{\bf Remark.}
It is clear from examples that in general this fibration is not a holomorphic map,
so we might think about resolving the indeterminacies by blowing-up $X' \rightarrow X$.
It would be interesting to see if this can be done in a way such that $X'$ has a split tangent bundle.

\begin{proposition}
\label{propositionrationalquotienttrivalintersection}
Let $X$ be a uniruled compact K\"ahler manifold with split tangent bundle $T_X =V_1 \oplus V_2$.
Let $Z$ be a general fibre of the rational quotient map, and suppose
that $T_Z \cap V_1|_Z$ or $T_Z \cap V_2|_Z$ is integrable.
If $T_Z \cap V_1|_Z = V_1|_Z$, 
the manifold $X$ 
has the structure of an analytic fibre bundle $X \rightarrow Y$
such that $T_{X/Y}=V_1$. 
\end{proposition}

{\bf Proof.} 
By theorem \ref{theoremfibration} the condition $T_Z \cap V_1|_Z = V_1|_Z$ implies that
there exists an almost holomorphic map
\merom{\phi_1}{X}{Y_1} such that the general fibre $F_1$ is rationally connected
and satisfies 
\[
T_{F_i} =  (T_{Z} \cap V_1|_Z)|_{F_i} = V_1|_{F_1}.
\]
It follows that $V_1$ is integrable and has a rationally connected leaf. We conclude with 
proposition \ref{propositionquotientsubmersion}. $\square$

In view of proposition \ref{propositionrationalquotienttrivalintersection} it is clear that theorem 
\ref{theoremcorank2} follows as soon as we have understood the geometry when $T_Z \cap V_1|_Z$ is a line bundle.
Since we will consider this situation also in the next section, we state this case as a

\begin{proposition}
\label{propositionuniruledsmallrank}
Let $X$ be a uniruled compact K\"ahler manifold 
with split tangent bundle $T_X = V_1 \oplus V_2$ where
$\rk V_1 = 2$. Suppose that the general fibre $Z$ of the rational quotient map 
satisfies $\rk (T_Z \cap V_1|_Z) = 1$.
Then there exists an equidimensional map \holom{\phi}{X}{Y} on a compact K\"ahler variety such that
the general fibre $F$ is a rational curve that satisfies $T_F \subset V_1|_F$.
\end{proposition}

{\bf Proof.} The line bundle $T_Z \cap V_1|_Z$ is integrable, so 
by theorem \ref{theoremfibration} there exists an almost holomorphic map
\merom{\phi_1}{X}{Y_1} such that the general fibre $F_1$ is rationally connected
and satisfies 
\[
T_{F_1} =  (T_{Z} \cap V_1|_Z)|_{F_1} \subset V_1|_{F_1}.
\]
Since $\rk (T_Z \cap V_1|_Z) = 1$, the general fibre is a smooth rational curve
such that
\[
T_X|_{F_1} \simeq \sO_{\PP^1} (2) \oplus \sO_{\PP^1}^{\dim X-1}.
\]
Since $h^0(F_1, N_{F_1/X})= \dim X -1$ and $h^1(F_1, N_{F_1/X})=0$,  
the corresponding open subvariety of the 
cycle space $\Chow{X}$ is smooth of dimension $\dim X-1$.
We denote by $Y$ its closure in $\Chow{X}$ and endow it with the reduced structure.
Denote by $\Gamma \subset Y \times X$ the reduction of the graph over $Y$. 
Denote furthermore by \holom{p_X}{\Gamma}{X} and \holom{p_Y}{\Gamma}{Y} the restrictions
of the projections to the graph.

{\it Step 1. We show that $p_X$ is finite.}
We argue by contradiction, then by the analytic version of Zariski's main theorem
there are fibres of positive dimension.
Let $x \in X$ be a point such that $\fibre{p_X}{x}$ has
a component of positive dimension. Let $\Delta \subset
p_Y(\fibre{p_X}{x})$ be an irreducible component of dimension $k>0$. 
Then $\Gamma_\Delta:=\fibre{p_Y}{\Delta}$ has dimension $k+1$. 
Consider now the foliation induced by $p_X^* V_1$ on $\Gamma \subset Y \times X$. Since 
a general
$p_{Y}$-fibre is contained in a $p_X^* V_1$-leaf and this is a closed condition, every fibre
$\fibre{p_Y}{y}$ is contained in a $p_X^* V_1$-leaf. 
So for $y \in \Delta$, the set $p_X(\fibre{p_Y}{y})$ 
is contained in $\mathfrak{V}_1^x$, the $V_1$-leaf through $x$.
It follows that $S:=p_X(\fibre{p_Y}{\Delta})$ is contained set-theoretically in $\mathfrak{V}_1^x$.
Since $p_X$ is injective on the fibres of $p_Y$, and \fibre{p_Y}{\Delta} has dimension $k+1 \geq 2$,
the subvariety $S$ has dimension at least 2.
Since $\rk V_1=2$, it has dimension 2 and $S=\mathfrak{V}_1^x$ (at least set-theoretically).
So $\mathfrak{V}_1^x$ is a compact leaf and is covered by a family of rational cycles that
intersect in the point $x$. Hence $\mathfrak{V}_1^x$ is rationally connected,
so by proposition \ref{propositionquotientsubmersion}
there exists a submersion \holom{\psi}{X}{Z}
such that $T_{X/Z}=V_1$ and the fibres are rationally connected.
By the universal property of the rational quotient the
general $\psi$-fibre is contracted by rational quotient map. 
This implies $\rk (T_Z \cap V_1|_Z) = \rk V_1$, a contradiction. 

{\it Step 2. Construction of $\phi$.} 
Since $p_X$ is birational by construction and finite, it is bijective 
by the analytic version of  
Zariski's main theorem. 
Since $X$ is
smooth and $\Gamma$ reduced this shows that $p_X$ is an isomorphism.  
Since $p_Y$ is equidimensional, \holom{\phi:=p_Y \circ p_X^{-1}}{X}{Y}
is equidimensional.
$\square$

\medskip

{\bf Proof of theorem \ref{theoremcorank2}.}
By proposition \ref{propositionrationallyconnectedungeneric}, the general fibre $Z$ of
the rational quotient map satisfies 
\[
T_Z = (V_1|_Z \cap T_Z) \oplus (V_2|_Z \cap T_Z).
\]
Since $\rk V_1=2$, there are three cases.

If $V_1|_Z \cap T_Z = V_1|_Z$, we conclude with proposition \ref{propositionrationalquotienttrivalintersection}.

If $0 \subsetneq V_1|_Z \cap T_Z \subsetneq T_Z$, the intersection has rank 1. Proposition 
\ref{propositionuniruledsmallrank} shows that we are in the second case of the statement.

If $V_1|_Z \cap T_Z=0$, clearly $T_Z = T_Z \cap V_2|_Z \subset V_2|_Z$. $\square$

\section{The projective case}

The main setback of theorem \ref{theoremcorank2} is that in the second case it is not clear
if the image of the constructed fibration is smooth.  
In order to refine our analysis of this fibration we have to use the theory of Mori contractions, 
this forces us leave the K\"ahler world. In lemma \ref{lemmauniruledintegrability} we will then show
the integrability of at least one direct factor for a splitting in vector bundles of rank 2. 
For uniruled varieties, the statement does not generalise to a splitting in vector bundles of higher rank.
Nevertheless the lemma provides some first evidence for conjecture \ref{conjecturerationallyconnectedintegrable}
which it establishes for manifolds of dimension at most four. 

A Mori contraction of a projective manifold $X$ is
a morphism with connected fibres \holom{\phi}{X}{Y} to a normal variety $Y$ such that the anticanonical bundle 
$-K_X$ is $\phi$-ample.
We say that the contraction is elementary if the relative Picard number
$\rho(X/Y)$ is equal to one.
The contraction is said to be of fibre type if $\dim Y < \dim X$; otherwise it is birational. 

\begin{lemma}
\label{lemmarelativedim1flattening}
Let $X$ be a projective manifold, and 
let \holom{\phi}{X}{Y} be an equidimensional fibration of relative dimension 1
on a normal variety $Y$
such that the general fibre $F$ is a rational curve.
Then there exists a factorisation $\phi=\tilde{\phi} \circ \mu$, where 
\holom{\mu}{X}{\tilde{X}} is a birational morphism onto a projective manifold $\tilde{X}$ and
\holom{\tilde{\phi}}{\tilde{X}}{Y} makes $\tilde{X}$ into a $\PP^1$- or conic bundle. 
Furthermore $\mu$ is a composition of blow-ups of projective manifolds along submanifolds
of codimension 2, and $Y$ is smooth.
\end{lemma}

{\bf Proof.}
We argue by induction on the relative Picard number $\rho(X/Y)$.
If $\rho(X/Y)=1$, the anticanonical divisor $-K_X$ is $\phi$-ample and the contraction is elementary,
so by Ando's theorem $\phi$ induces a $\PP^1$-bundle or a conic bundle structure. In both cases $Y$ is smooth.

Suppose now that $\rho(X/Y)>1$. Since the general fibre is a rational curve, 
the canonical divisor is not $\phi$-nef.
It follows from the relative contraction theorem \cite[Thm.4-1-1]{KMM87} that there exists an elementary contraction 
\holom{\mu}{X}{X'} that is a $Y$-morphism, i.e. there exists a morphism \holom{\tilde{\phi}}{\tilde{X}}{Y}
such that $\phi = \tilde{\phi} \circ \mu$. Since $\phi$ is equidimensional of relative dimension 1, it follows
that all the $\mu$-fibres have dimension at most 1. Thus $\mu$ is of fibre type of relative dimension 1 
or of birational type. 

We claim that $\mu$ is not of fibre type.
We argue by contradiction and suppose that $\dim X = \dim \tilde{X}+1$. 
Then $\dim \tilde{X}=\dim Y$, so $\tilde{\phi}$ is a birational morphism. Since
$\rho(\tilde{X}/Y)=\rho(X/Y)-1>0$, the map $\tilde{\phi}$ is not an isomorphism, 
so there exists a fibre $\fibre{\tilde{\phi}}{y}$
of positive dimension. Since $\mu$ is of fibre type, we see that
$\fibre{\phi}{y}= \fibre{\mu}{\fibre{\tilde{\phi}}{y}}$ has dimension at least 2, a contradiction.    

Hence $\mu$ is a birational contraction such that all the fibres have dimension at most 1.
A straightforward application of  Wi\'sniewski's inequality 
\cite[Thm. 1.1]{Wi91} shows that the contraction is divisorial 
and all the fibres have dimension at most 1.
By Ando's theorem \cite[Thm. 2.1]{An85}
we know that $\tilde{X}$ is smooth and $\mu$ is the blow-up of $\tilde{X}$ 
along a smooth submanifold of codimension 2.
Now $\rho(\tilde{X}/Y)=\rho(X/Y)-1$ and $\tilde{\phi}$ is equidimensional of relative dimension 1 over $Y$, so 
the statement follows by the induction hypothesis.
$\square$

\medskip

{\bf Remark.} 
In order to generalise the proof to the compact K\"ahler case it would be necessary to establish
a relative contraction theorem for projective morphisms between compact K\"ahler varieties.
Unfortunately the Mori theory for compact K\"ahler manifolds is not yet at this stage, in particular
there seem to be no statements for the relative situation.

\begin{corollary}
\label{corollaryprojectivecorank2} 
Let $X$ be a uniruled projective manifold such that $T_X = V_1 \oplus V_2$ and $\rk V_1=2$.
Let $Z$ be a general fibre of the rational quotient map, and suppose
that $T_Z \cap V_1|_Z$ or $T_Z \cap V_2|_Z$ is integrable.
Suppose that $T_Z \cap V_1|_Z \neq 0$. Then $X$ admits a flat fibration $\holom{\phi}{X}{Y}$
on a smooth projective manifold $Y$ such that 
\[
T_Y= (\phi_* (T \phi (V_1)))^{**} \oplus (\phi_* (T \phi (V_2)))^{**}.
\]
\end{corollary}
 
{\bf Proof.}
If $T_Z \cap V_1|_Z=V_1|_Z$ we conclude with the first case of theorem \ref{theoremcorank2}. 

If $0 \subsetneq T_Z \cap V_1|_Z \subsetneq T_Z$ we are in the second case of theorem \ref{theoremcorank2},
so there exists an equidimensional fibration $X \rightarrow Y$ such that the general fibre is a rational curve.
Since $X$ is projective, there exists by lemma \ref{lemmarelativedim1flattening} 
a factorisation $\phi=\tilde{\phi} \circ \mu$, where 
\holom{\mu}{X}{\tilde{X}} is birational morphism onto a projective manifold $\tilde{X}$ and
\holom{\tilde{\phi}}{\tilde{X}}{Y} makes $\tilde{X}$ into a $\PP^1$- or conic bundle. 
Since $\mu$ is a composition of smooth blow-ups a repeated application of lemma \ref{lemmasplittingblowup}
shows that
\[
T_{\tilde{X}} = (\mu_* V_1)^{**} \oplus (\mu_* V_2)^{**},
\]
so we can apply lemma \ref{lemmapushdown} to $\tilde{\phi}$ to see that for $i=1,2$
\[
W_i := (\phi_* (T \phi (V_i)))^{**} = (\tilde{\phi}_* (T \tilde{\phi} ((\mu_* V_i)^{**})))^{**} 
\]
is a subbundle of $T_Y$ such that $T_Y=W_1 \oplus W_2$.
$\square$

\begin{lemma}
\label{lemmauniruledintegrability}
Let $X$ be a uniruled projective manifold such that $T_X= \oplus_{j=1}^k V_j$, where
for all $j=1, \ldots, k$ we have $\rk V_j \leq 2$. Then one of the direct factors is integrable.

In particular if $\dim X \leq 4$, one of the direct factors is integrable.
\end{lemma}

\medskip

{\bf Proof.}
The statement is trivial if one direct factor has rank 1, so we suppose that all the direct factors have rank 2.
Let \holom{f}{\PP^1}{X} be a general minimal rational curve on $X$, then 
\[
\bigoplus_{j=1}^k f^* V_j = f^* T_X \simeq \sO_{\PP^1}(2) \oplus \sO_{\PP^1}(1)^{\oplus a} 
\oplus \sO_{\PP^1}^{\oplus b},
\]
We may suppose up to renumbering that $f^* V_1 \simeq \sO_{\PP^1}(2) \oplus \sO_{\PP^1}(c)$ where $c=0$ or $1$.
It follows that for $i \geq 2$, we have $f^* V_i \simeq  \sO_{\PP^1}(1) \oplus \sO_{\PP^1}$
or $f^* V_i \simeq  \sO_{\PP^1}(1)^{\oplus 2}$ or  $f^* V_i \simeq  \sO_{\PP^1}^{\oplus 2}$,
in particular
\[
H^1( \PP^1, f^* V_i^*) = 0 \qquad \forall \ i \geq 2.
\]
By \cite[Lemma 0.4]{CP02}, we have $c_1(V_i) \in H^1(X,V_i^*)$, so $c_1(f^*V_i) \in H^1(\PP^1, f^* V_i^*)$ 
is zero for $i \geq 2$. So $f^* \det V_i \simeq \sO_{\PP^1}$, since $f^* V_i$ is nef this implies
$f^* V_i \simeq  \sO_{\PP^1}^{\oplus 2}$ for $i \geq 2$.
It follows that $\bigwedge^2 V_1|_{f(\PP^1)}$ is ample and $(T_X/V_1)|_{f(\PP^1)}= \oplus_{j \geq 2} V_j|_{f(\PP^1)}$ 
is trivial. By lemma \ref{lemmaintegrability} this implies the integrability of $V_1$. $\square$

\medskip

{\bf Proof of theorem \ref{theoremrank2}.}
Let $Z$ be a general fibre of the rational quotient map, then
an application of proposition \ref{propositionrationallyconnectedungeneric}
to all the possible decompositions $W_1:=V_i$ and $W_2:= \oplus_{j=1, j \neq i}^k V_j$
implies
\[
T_Z = \oplus_{j=1}^k (T_Z \cap V_j|_Z).
\]
Furthermore $\rk (T_Z \cap V_j|_Z) \leq 2$ for all $j \in \{1, \ldots, k\}$, so
one of the direct factors is integrable by lemma \ref{lemmauniruledintegrability}.
We can  now apply theorem \ref{theorembogomolovbetter} inductively to see that all the direct factors
are integrable and $Z$ splits in a product. 

{\it Step 1. Integrability of the direct factors.}
Suppose that $T_Z \cap V_i|_Z$ is
not zero, then there are two possibilities. Either $T_Z \cap V_i|_Z = V_i|_Z$, so the integrability
of $V_i$ follows from the integrability of $V_i|_Z$. Or $T_Z \cap V_i|_Z \subsetneq V_i|_Z$, 
then $V_i$ has rank 2 and the splitting of $Z$ in a product yields a dominant family of rational curves such 
that a general member $C$ satisfies $T_C \subset V_i|_C$ and the normal bundle $N_{C/X}$ is trivial.
Since 
\[
T_X|_C = V_i|_C \oplus \bigoplus_{j=1, j \neq i}^k V_j|_C = T_C \oplus N_{C/X}
\]
and $\rk V_i=2$, this implies that $(\wedge^2 V_i)|_C$ is ample and 
$(T_X/V_i)|_C \simeq \bigoplus_{j=1, j \neq i}^k V_j|_C$ is trivial.
By lemma \ref{lemmaintegrability}
this implies the integrability of $V_i$.

{\it Step 2. Structure of the rational quotient map.} 
We proceed
by induction on the dimension of $X$, the case $\dim X=1$ is trivial.
Up to renumbering we can suppose that the intersection $T_Z \cap V_1|_Z$ is not empty. 

If $V_1$ has rank 1, we have $T_Z \cap V_1|_Z=V_1|_Z$, so $V_1$ is integrable and the general leaf
is rationally connected. Thus by proposition \ref{propositionquotientsubmersion} there exists a submersion
\holom{\psi}{X}{Y'} such that $T_{X/Y'}=V_1$. Hence $T_Y'= \oplus_{j=2}^k (\psi_* (T \psi (V_j)))^{**}$. 
If $Y'$ is not uniruled we are done, otherwise apply the induction
hypothesis to $Y'$. 

If $V_2$ has rank 2, we apply corollary \ref{corollaryprojectivecorank2} to obtain a flat
fibration \holom{\psi}{X}{Y'} onto a projective manifold $Y'$ 
such that the general fibre is rationally connected and
\[
T_{Y'} = 
(\psi_* (T \psi (V_1)))^{**} 
\oplus 
(\psi_* (T \psi (\oplus_{j=2}^k V_j)))^{**} 
= \oplus_{j=1}^k 
(\psi_* (T \psi (V_j)))^{**} 
\]
If $Y'$ is not uniruled we are done, otherwise apply  the induction hypothesis to $Y'$.

{\it Step 3. Integrability of the images.}
Let \holom{\phi}{X}{Y} be the map constructed in step 2. 
Then
\[
T_{Y} = \oplus_{j=1}^k (\phi_* (T \phi (V_j)))^{**} 
\]
and $Y$ is not uniruled. Therefore for all $j \in \{1, \ldots, k\}$, the vector bundle
$(\phi_* (T \phi (V_i)))^{**}$ is integrable by \cite[Thm. 1.3]{H05}.
$\square$

\section{An application to the universal covering}
\label{sectionuniversalcovering}

This section is essentially devoted to the proof 
of theorem \ref{theoremcorank2projective}. The basic strategy
is to prove conjecture \ref{conjecturebeauville} by a reduction to the case of non-uniruled varieties 
and induction on the rank of the direct factors.
Before we come to the proof we have to show a refinement of \cite[Thm.1.3]{H05}.

\begin{lemma} 
\label{lemmanonuniruledintegrability}
Let $X$ be a projective manifold with split tangent bundle $T_X=V_1 \oplus V_2$.
Suppose that a general fibre $Z$ of the rational quotient map satisfies $T_Z \subset V_2|_Z$.
Then $V_2$ is integrable and $\det V_1^*$ is pseudo-effective.
\end{lemma} 

{\bf Proof.} 
{\it Step 1.}
Suppose that $L:=\det V_1^*$ is pseudoeffective.
Since $V_1^*$ is a direct factor of $\Omega_X$, the vector bundle $\det V_1 \otimes \wedge^{\rk V_1} \Omega_X$
has a trivial direct factor. 
If $\theta \in H^0(X, L^{-1} \otimes \wedge^{\rk V_1} \Omega_X)$ is the associated nowhere-vanishing
$\det V_1$ -valued form, and $\zeta$ a germ of any vector field, a local computation shows that $i_\zeta \theta = 0$
if and only if $\zeta$ is in $V_2$. 
An integrability criterion by Demailly \cite[Thm.]{De02} shows that $V_2$ is integrable.

{\it Step 2. $\det V_1^*$ is pseudoeffective.}
We argue by contradiction, then by \cite{BDPP04} there exists a birational morphism
\holom{\phi}{X'}{X} and a general intersection curve $C:=D_1 \cap \ldots \cap D_{\dim X-1}$  of very ample divisors 
$D_1, \ldots, D_{\dim X-1}$ where $D_i \in |m_i H|$ for some ample divisor $H$
such that $\phi^* \det V_1^* \cdot L < 0$. 
Let
\[
0 = E_0 \subset E_1 \subset \ldots \subset E_r= \phi^* V_1
\]
be the Harder-Narasimham filtration with respect to the polarisation $H$, i.e. the graded pieces
$E_{i+1}/E_i$ are semistable with respect to $H$. 
Since $m_1, \ldots, m_{\dim X-1}$ can be arbitrarily high,
we can suppose that the filtration commutes with restriction to $C$. 
Furthermore since $C$ is general and $E_1$ a reflexive sheaf, the curve $C$ is contained
in the locus where $E_1$ is locally free. Since 
\[
\mu(E_1|_C) \geq \mu(\phi^* V_1|_C) =  \frac{\deg_C V_1}{\rk \phi^* V_1} > 0
\]
and $E_1|_C$ is semistable, it is ample by \cite[p.62]{Laz04a}.
By \cite[Cor.1.5]{KST05} 
this implies that $E_1$ is vertical with respect to the rational quotient map,
that is a general fibre $Z$ of the rational quotient satisfies $E_1|_Z \cap T_Z = E_1|_Z$.
It follows that the intersection $T_Z \cap V_1|_Z$ is not zero, a contradiction. 
$\square$

\medskip

{\bf Proof of theorem \ref{theoremcorank2projective}.}
If $T_Z \cap V_1|_Z=0$, proposition \ref{propositionrationallyconnectedungeneric} 
shows that $T_Z \subset V_2|_Z$. Therefore we can conclude with lemma \ref{lemmanonuniruledintegrability}.

If $T_Z \cap V_1|_Z=V_1|_Z$ there exists a submersion $X \rightarrow Y$ such that $T_{X/Y}=V_1$.
Furthermore $V_2$ is an integrable connection on the submersion, so we conclude with the Ehresmann
theorem \cite[Thm. 3.17]{H05}.  

If $T_Z \cap V_1|_Z \subsetneq V_1|_Z$ is a line bundle there exists 
an equidimensional map $\holom{\phi}{X}{Y}$ of relative dimension one
such that the general $\phi$-fibre $F$ is a rational curve and $T_F \subset V_1|_F$. 
Since $X$ is projective there exists by lemma \ref{lemmarelativedim1flattening} a 
a factorisation $\phi=\tilde{\phi} \circ \mu$, where 
\holom{\mu}{X}{\tilde{X}} is birational morphism onto a projective manifold $\tilde{X}$ and
\holom{\tilde{\phi}}{\tilde{X}}{Y} makes $\tilde{X}$ into a $\PP^1$- or conic bundle. 
Furthermore $\mu$ is a composition of blow-ups of projective manifolds along submanifolds
of codimension 2, so lemma \ref{lemmasplittingblowup} implies
\[
T_{\tilde{X}} = (\mu_* V_1)^{**} \oplus (\mu_* V_2)^{**}.
\]
By the same lemma 
it is sufficient to show the conjecture for $\tilde{X}$, so we can replace without loss of generality
$X$ by $\tilde{X}$ and suppose that the fibration $\phi$ makes $X$ into a $\PP^1$- or conic bundle
over the projective manifold $Y$. Set $W_j:= (\phi_* (T \phi(V_j)))^{**}$, then  
\[
T_Y = W_1 \oplus W_2
\]
by lemma \ref{lemmapushdown} and $W_1$ has rank 1.
Furthermore by (\cite[Prop. 4.23.]{H05}, see also \cite[Cor. 4.3.9]{H06}) all the fibres of $\phi$ are reduced.

The manifold $Y$ can't have the structure of a $\PP^1$-bundle $Y \rightarrow M$ such that 
$T_{Y/M}=W_1$: this would yield a morphism $X \rightarrow M$ such that $T_{X/M}=V_1$
and the general fibre is rationally connected. This contradicts $T_Z \cap V_1|_Z \subsetneq V_1|_Z$.
Therefore by \cite[Thm. 1]{BPT04}
the subbundle $W_2$ is integrable, 
and the universal covering \holom{\mu}{\tilde{Y}}{Y} satisfies $\tilde{Y} \simeq Y_1 \times Y_2$
such that $\mu^* W_1 = p_{Y_1}^* T_{Y_1}$  and $\mu^* W_2 = p_{Y_2}^* T_{Y_2}$. 
Furthermore we have a commutative diagram
\[
\xymatrix{
\tilde{X} \ar@/_/[dd]_q \ar[r]^{\tilde{\mu}} \ar[d]^{\tilde{\phi}} & X  \ar[d]^\phi
\\
\tilde{Y}  \ar[d]^{p_{Y_2}} \ar[r]^\mu & Y
\\
Y_2
}
\]
where \holom{\tilde{\mu}}{\tilde{X}:= X \times_Y \tilde{Y}}{X} is \'etale.
By construction  the set-theoretical fibres of $q$
are $\tilde{\mu}^* V_1$-leaves. Since $\phi$ has no multiple fibres, the fibration
$\tilde{\phi}$ has no multiple fibres. Hence $q=p_{Y_2} \circ \tilde{\phi}$ does not have
any multiple fibres, so the scheme-theoretical fibres are $\tilde{\mu}^* V_1$-leaves. This shows
that $q$ is a submersion with integrable connection $\tilde{\mu}^* V_2$.
Since
\[
(\tilde{\phi}_* (T \tilde{\phi} \tilde{\mu}^* V_2))^{**} = \mu^* W_2 = p_{Y_2}^* T_{Y_2},
\]
there exists for every $\tilde{\mu}^* V_2$ leaf $\mathfrak{V}_2$ a $y_1 \in Y_1$ such that
$\tilde{\phi}(\mathfrak{V}_2)= y_1 \times Y_2$. By lemma \ref{lemmaehresmann} below the restriction
of $q$ to a   $\tilde{\mu}^* V_2$ leaf is an \'etale covering, so we conclude with the Ehresmann
theorem \cite[Thm. 3.17]{H05}. $\square$

{\bf Remark.}
Note that theorem \ref{theoremcorank2projective} generalises immediately to the compact K\"ahler case if we show that
the map in the second case of theorem \ref{theoremcorank2} is flat.
 
\begin{lemma}
\label{lemmaehresmann}
Let \holom{\phi}{X}{Y_1 \times Y_2} be a proper surjective map  
from a complex manifold $X$ onto a product of (not necessarily compact) complex manifolds
such that the morphism \holom{q:=p_{Y_2} \circ \phi}{X}{Y_2} is a submersion that admits
an integrable connection $V \subset T_X$. Suppose that for every $V$-leaf $\mathfrak{V}$, there
exists a $y_1 \in Y_1$ such that $\phi(\mathfrak{V})= y_1 \times Y_2$. 
Then the restriction of $q$ to every $V$-leaf
is an \'etale covering.
\end{lemma}

The proof consists merely of rephrasing the classical proof of the Ehresmann theorem
as in \cite[V.,\S 2,Prop.1]{CLN85}. For the convenience of the reader we nevertheless include this
technical exercise.

{\bf Proof.} In this proof all fibres and intersections are set-theoretical.

Let $\mathfrak{V}$ be a $V$-leaf, and let $y_1 \in Y_1$ such that $\phi(\mathfrak{V})=y_1 \times Y_2$.
Since \holom{p_{Y_2}|_{y_1 \times Y_2}}{y_1 \times Y_2}{Y_2} is an isomorphism, it is sufficient
to show that \holom{\phi|_\mathfrak{V}}{\mathfrak{V}}{y_1 \times Y_2} is an \'etale map. 
Furthermore it is sufficient to show that for $y_1 \times y_2 \in y_1 \times Y_2$, 
there exists a disc $D \subset y_1 \times  Y_2$
such that for $y \in D$, the fibre \fibre{\phi}{y} cuts each leaf of the restricted
foliation $V|_{\fibre{\phi}{D}}$ exactly in one point. Granting this for the moment, 
we show how this implies the result.  
The connected components of $\mathfrak{V} \cap \fibre{\phi}{D}$ are leaves of $V|_{\fibre{\phi}{D}}$ 
Let $\mathfrak{V}'$ be such a connected component. Since for $y \in D$, the intersection
$\mathfrak{V}' \cap \fibre{\phi}{y}$ is exactly one point, the restricted morphism
\holom{\phi|_{\mathfrak{V}' }}{\mathfrak{V}' }{D} is one-to-one and onto, so it is a biholomorphism.
This shows that  \holom{\phi|_{\mathfrak{V} \cap \fibre{\phi}{D}}}{\mathfrak{V} \cap \fibre{\phi}{D}}{D} 
is a trivialisation of $\phi|_{\mathfrak{V}}$.

Let us now show the claim. Set $k:= \rk V$ and $n:=\dim X$, and set $Z:=\fibre{\phi}{y_1 \times Y_2}$. 
Since every $V$-leaf is sent on some $b \times Y_2$, the complex space $Z$ is $V$-saturated.
In particular if $\mathfrak{V} \subset Z$ is leaf,
the restriction of a distinguished map \holom{f_i}{W_i}{\D^{n-k}}
to $Z$ which we denote by \holom{f_i|_{W_i \cap Z}}{W_i \cap Z}{\D^{n-k}},
is a distinguished map for the foliation $V|_Z$ and a plaque of $f_i$ is contained in $\mathfrak{V}$
if and only if it is a plaque of $f_i|_{W_i \cap Z}$.

{\it Step 1. The local situation.} Let $x \in \fibre{\phi}{y_1 \times y_2}$ be a point.
Since $q$ is a submersion with integrable connection $V$ there exists coordinate 
neighbourhood $x \in W_x' \subset X$ with local coordinates $z_1, \ldots, z_{k}, z_{k+1}, \ldots, z_{n}$
and a coordinate neighbourhood $y_2 \in U_{x} \subset Y_2$ with coordinate $w_1, \ldots, w_k$
such that $q(W_x')=U_x$ and \holom{q|_{W_x'}}{W_x}{U_{x}} is given in these coordinates
by 
\[
(z_1, \ldots, z_n) \rightarrow (z_1, \ldots, z_k).
\] 
Furthermore there exists a distinguished
map \holom{f_x}{W_x'}{\D^{n-k}} given in these coordinates by
\[
(z_1, \ldots, z_n) \rightarrow (z_{k+1}, \ldots, z_n).
\]
Since $x \in \fibre{\phi}{y_1 \times y_2}$ and $\phi$ is equidimensional over a smooth base,
so open, $\phi(W_x')$ is a neighbourhood of $y_1 \times y_2$ in $Y_1 \times Y_2$.
Since \holom{p_{Y_2}|_{y_1 \times U_x}}{y_1 \times U_x}{U_x} is an isomorphism 
we can suppose that up to restricting $U_x$ and $W_x'$ a bit that 
\[
\phi(W_x') \cap (y_1 \times Y_2) = y_1 \times U_x.
\]
Set $W_x:= W_x' \cap Z$, then $\phi|_Z(W_x)= y_1 \times U_x$.
It then follows from this local description 
that \holom{\phi|_{W_x}}{W_x}{y_1 \times U_x} has the property that for
$y \in y_1 \times U_x$ the fibre \fibre{\phi}{y} intersects each plaque of the distinguished map 
$\holom{f_x|_{W_x}}{W_x}{\D^k}$
in exactly one point. 

{\it Step 2. Using the properness.} 
Since the fibre \fibre{\phi}{y_1 \times y_2} is compact,
we can take a finite cover of the fibre by $W_i:=W_{x_i}$ where $i=1, \ldots, l$ 
and $W_{x_i}$ is as in step 1. 
For each $i \in \{ 1, \ldots, l\}$, the image $\phi(W_i)$ is a neighbourhood of 
$y_1 \times y_2 \in y_1 \times Y_2$. Let $D \subset \cap_{i=1}^l \phi(W_i)$
be a disc that contains $y_1 \times y_2$.
If $\mathfrak{V}'$ is a leaf of $V|_{\fibre{\phi}{D}}$, it is contained
in some plaque $P$ of $W_i$ for some $i \in \{ 1, \ldots, l\}$.
Since the plaques intersect each fibre at most in one point, 
\holom{\phi|_{\mathfrak{V}'}}{\mathfrak{V}'}{D} is injective. The equality 
$P \cap \fibre{\phi}{D} = \mathfrak{V}'$ then implies that
\[
\phi (\mathfrak{V}') = \phi (P \cap \fibre{\phi}{D}) = \phi(P) \cap D = D, 
\] 
so \holom{\phi|_{\mathfrak{V}'}}{\mathfrak{V}'}{D} is surjective. 
So $\mathfrak{V}'$ intersects each fibre exactly in one point.
$\square$

\end{document}